\documentclass[final,a4paper,notitlepage,12pt,fullpage]{article}
\usepackage[dvips]{graphics}
\usepackage{latexsym}
\usepackage{amsfonts}
\usepackage{amssymb}
\usepackage{amsmath}
\usepackage[dvips]{graphics}
\usepackage{theorem}
\usepackage{epsfig}

\newenvironment{@abssec}[1]{%
     \if@twocolumn
       \section*{#1}%
     \else
       \vspace{.05in}\footnotesize
       \parindent .2in
         {\bfseries #1. }\ignorespaces
     \fi}
     {\if@twocolumn\else\par\vspace{.1in}\fi}
\newenvironment{keywords}{\begin{@abssec}{Key words}}{\end{@abssec}}
\newenvironment{AMS}{\begin{@abssec}{AMS subject classification}}{\end{@abssec}}

\def\eqbd{\mathop{{:}{=}}}

\def\esssup{{\hbox{\rm esssup}}}
\newcommand{\R}{\mathbb R}
\newcommand{\C}{\mathbb C}
\newcommand{\F}{\mathbb F}
\newcommand{\N}{\mathbb N}

\newcommand{\inn}{\! \in \!}

\newtheorem{theo}{Theorem}[section]
\newtheorem{cor}[theo]{Corollary}

\newtheorem{lemma}[theo]{Lemma}
\newtheorem{prop}[theo]{Proposition}
\newenvironment{proof}{ \par \noindent {\bf Proof:} \normalsize}
{\hfill {$\square$} \break \par\vspace{.3cm}}
\newenvironment{proofof}{ \par \noindent {\bf Proof of} \normalsize}
{\hfill {$\square$} \break \par\vspace{.3cm}}
\begin{document}

%
%
\title{Real and Complex Operator Norms}
\author{ Olga Holtz\thanks{Supported 
by the DFG center ``Mathematics for key technologies'' in Berlin, Germany.} \\
University of California-Berkeley\\
{\tt holtz@Math.Berkeley.EDU}\\[.3 cm]
Michael Karow  \\
Technische Universit\"at Berlin \\
{\tt karow@math.TU-Berlin.de}
\\}
\date{July 2004}
\maketitle

%
%



\begin{abstract} Real and complex norms of a linear operator 
acting on a normed complexified space are considered. Bounds
on the ratio of these norms are given. The real and complex 
norms are shown to coincide for four classes of operators:
\begin{enumerate}
\item real linear operators from $L_p(\mu_1)$ to $L_q(\mu_2)$, 
$1\leq p\leq q\leq \infty$;
\item real linear operators between inner product spaces;
\item nonnegative linear operators acting between complexified 
function spaces with absolute and monotonic norms;
\item real linear operators from a complexified function space
with a norm satisfying $\|\Re x \|\leq \|x\|$  to $L_\infty(\mu)$. 
\end{enumerate}
The inequality $p\leq q$ in Case 1 is shown to be sharp.

A class of norm extensions from a real vector space to
its complexification is constructed that preserve operator norms. 
\end{abstract}

%
%
\section{Introduction}
By a normed complexified vector space we mean a  triple $(X,X_\R,\|\cdot\|_X)$, where 
\begin{itemize}
\item $X$ and $X_\R$ are  vector spaces over $\C$ and $\R$, respectively;
\item  $X$ is the algebraic complexification of $X_\R$, i.e. each 
$x \inn X$ can be uniquely written in the form $x=x_1+i x_2$ with $x_1,x_2 \inn X_\R$, we set $x_1\eqbd \Re x$, $x_2\eqbd \Im x$;
\item the function $\|\cdot\|_X:X\to [0, \infty)$ is a norm on the complex vector space $X$, in particular $\|\lambda \,x \|_X=|\lambda|\, \|x\|_X$ for all $\lambda \inn \C$, $x\inn X$.  
\end{itemize}
Let $(X,X_\R,\|\cdot\|_X)$  be a normed complexified space, let  $(Y,\|\cdot\|_Y)$
be a normed vector space over $\C$,
 and  $A:X \to Y$ be a $\C$-linear operator.  
Then two operator norms of $A$ with respect to
the norms $\|\cdot\|_X$, $\|\cdot\|_Y$ can be defined
  depending on whether the
 supremum of the quotient $\|Ax\|_Y/\|x\|_X$ is taken over all 
$x \inn X\setminus \{0\}$ or  over all $x \inn X_\R \setminus \{0\}$.
We distinguish these norms by a superscript $\C$ and $\R$ respectively,
 i.e. we set
$$\|A\|_{X, Y}^{\C}\eqbd
\sup_{x \in X\setminus\{0\}}\frac{\|Ax\|_Y}{\|x\|_X},
\qquad
\|A\|_{X, Y}^{\R}\eqbd
\sup_{x \in X_\R\setminus\{0\}}\frac{\|Ax\|_Y}{\|x\|_X}.$$
The case that the suprema are infinity is not excluded.
Obviously, $\|A\|_{X, Y}^{\R}\leq \|A\|_{X, Y}^{\C}$. 

In this note we address the question how large the ratio of these operator norms 
can be and under which conditions they coincide. Our main interest is in
the case that  $Y$ is the algebraic complexification of a real vector space $Y_\R$
and  $A$ is a {\em real operator,\/} i.e.,
$A(X_\R) \subseteq Y_\R$.  We discuss the above-mentioned question 
in detail for inner product spaces and the Banach spaces 
$L_p^\F(\mu)=L_p^\F(\Omega, {\cal B}, \mu)$, where 
$p\inn [1, \infty]$, $\mu$ is a positive measure on a $\sigma$-algebra ${\cal B}$ of subsets of the set $\Omega$ and
$$L_p^\F(\mu)=\{\; x:\Omega \to \F \quad |  \quad x\text{ is $\cal B$-measurable}, \;\; \|x\|_p < \infty\,\}, \qquad \F=\R\text{ or }\C,$$
\begin{equation*}
\|x\|_p   \eqbd \begin{cases} \left( \int_\Omega |x(t)|^p\, d\mu_t\;\right)^{1/p}, & 1 \leq p <\infty,  \cr
\esssup_{t \in \Omega} \,|x(t)|, & p=\infty. \end{cases}
\end{equation*}
Recall that this definition includes as special cases the finite- and 
infinite-dimensional  $l_p$ spaces. They correspond to $\Omega=\{1, \ldots ,n\}$ 
or $\Omega=\N$ and the counting measure $\mu$ satisfying $\mu(\{t\})=1$ for all
$t \inn \Omega$. The associated norms are 
$$ \|x\|_p =\begin{cases} \left(\sum_{t\in \Omega}^n|x_t|^p\right)^{1/p},
&  1\leq p < \infty,\cr 
\sup_{t\in \Omega}|x_t|, & p=\infty, \end{cases} $$
where $x=[x_1 \ldots x_n]' \inn \C^n$ or $x=(x_j)_{j \in \N} \in l_p$. 

Bounds on the ratio  $\|A\|_{X,Y}^{\C}/\|A\|_{X,Y}^{\R}$ are discussed 
in Section~\ref{sec:bounds}. We show that it is bounded by $2$ 
whenever the norm on the source space satisfies the condition
$ \| \Re x\| \leq \|x\|.$ Tighter bounds are derived for a linear map 
acting on $L_p (\mu)$. If both the source and the target are 
inner product spaces, the ratio is shown to be bounded by $\sqrt{2}$.

Section~\ref{sec:equal} is devoted to classes of operators for which the 
real and the complex norms coincide. The main result of that section,
and of the paper, Theorem~\ref{theo:main}, shows that the norms are
equal for all real linear maps from $L_p(\mu_1)$ to $L_q(\mu_2)$,
provided that $1\leq p\leq q\leq \infty$. We discuss three more classes 
of linear operators whose real and complex norms are equal at the end
of Section~\ref{sec:equal}.

The goal of Section~\ref{sec:unequal} is to show that the main result of 
the paper cannot be improved, that is, to describe counterexamples to
Theorem~\ref{theo:main} in the case $p<q$ ($1\leq p, q\leq \infty$).

Perhaps the simplest example that shows that the real and complex norms 
of a real linear map may be different can be found already in $\C^2$.
Let $A=\left[ \begin{array}{rr} 1 & -1 \\ 1 & 1 \end{array} \right]$. 
Then, for every $x=[x_1\;x_2]'\inn \C^2 \setminus \{0\}$,
$$\text{ }\qquad \frac{\|Ax\|_1}{\|x\|_\infty}=
\frac{|x_1-x_2|+|x_1+x_2|}{\max\{|x_1|,|x_2|\}}.$$
It can be checked by straightforward calculation~\cite{Book} that
$\|A\|_{\infty, 1}^\R=2$ but $\|A\|_{\infty, 1}^\C= |1-i|+|1+i|=2\sqrt{2}$.
It turns out that a finite-dimensional variation on this example covers the 
whole range $p<q$, as discussed in detail in Section~\ref{sec:unequal}.

Finally, in Section~\ref{sec:extend} we construct a class of norm extensions 
from a real space to its complexification that preserve operator norm.

\section{Bounds for $\|A\|_{X,Y}^{\C}/\|A\|_{X,Y}^{\R}$}\label{sec:bounds}
In this section we give simple bounds for the ratio
$\|A\|_{X,Y}^{\C}/\|A\|_{X,Y}^{\R}$.
\begin{prop}\label{prop:cx}
Let $(X,X_\R,\|\cdot\|_X)$  be a normed complexified space. Then for any normed space $(Y,\|\cdot\|_Y)$  
and any linear operator $A:X \to Y$,
\begin{eqnarray}\label{eq:quotient}
\|A\|_{X,Y}^{\C}&\leq& c_X\, \|A\|_{X, Y}^{\R},
\end{eqnarray}
where
$$c_X\eqbd \sup_{0\not = x \in X}
\frac{\|\Re x\|_X+\|\Im x\|_X}{\|x\|_X}.$$
\begin{itemize}
\item[$(a)$] Suppose $\|\cdot\|_X$ satisfies 
\begin{eqnarray}\label{eq:Re_condition}
\|\Re x\|_X \leq \|x\|_X\qquad \text{\rm for all } \; x \inn X.
\end{eqnarray}
Then $c_X\leq 2$.
\item[$(b)$] Suppose the space $X_\R$ is endowed with an inner product 
$\langle \cdot, \cdot \rangle\!:\! X_\R \times X_\R\to \R$. Let $\|\cdot\|_X$
be the norm induced by the complexification of this inner product,
i.e. ,
\begin{eqnarray}\label{eq:inner-prod-norm}
\|x\|_X^2=\langle \Re x ,\Re x\rangle +  \langle \Im x ,\Im x\rangle
\qquad \text{\rm for all }\; x \inn X.
\end{eqnarray}
Then $c_X=\sqrt{2}$.
\end{itemize}
\end{prop}
\begin{proof}
For all $x \inn X$, we have
\begin{eqnarray*}
\|Ax\|_Y &=&\|A\, \Re x + i \, A\, \Im x\|_Y\\
& \leq &
\|A\, \Re x\|_Y+\|i\, A \, \Im x \|_Y \\
&=& \|A\, \Re x\|_Y+\|A\, \Im x \|_Y\\
&\leq & 
\|A\|_{X,Y}^{\R}\, (\|\Re x\|_X +\|\Im x\|_X).
\end{eqnarray*}
Thus
$$\frac{\|Ax\|_Y}{\|x\|_X }\leq
\|A\|_{X,Y}^{\R}\frac{\|\Re x\|_X +\|\Im x\|_X}{\|x\|_X}
\leq \|A\|_{X,Y}^{\R}\, c_X .$$
This gives inequality (\ref{eq:quotient}).\\
$(a).$ Condition (\ref{eq:Re_condition}) yields
$$\|\Re x\|_X+\|\Im x\|_X=\|\Re x\|_X+\|\Re (ix)\|_X\leq \|x\|_X+\|ix\|_X=2\|x\|_X.$$ Thus $c_X \leq 2$.\\
$(b)$. Relation (\ref{eq:inner-prod-norm}) implies that
$\|x\|_X^2=\|\Re x \|_X^2+\|\Im x\|_X^2\geq (1/2)(\|\Re x \|_X+\|\Im x\|_X)^2$
 for all $x \inn X$. Thus $c_X\leq \sqrt{2}$. 
Let $x=(1+i)x_0$ with $x_0 \inn 
X_\R$. Then $\|x\|_X=
\sqrt{2}(\|\Re x \|_X+\|\Im x\|_X)$.
Thus $c_X \geq \sqrt{2}$.
 \end{proof}
A norm $\|\cdot\|_X$ is said to be conjugation-invariant if 
$\|x\|_X=\|\overline{x}\|_X$ for all $x\inn X$, where $\overline{x}
\eqbd\Re x -i \Im x$.
It is easily seen that a conjugation-invariant norm satisfies  
condition~(\ref{eq:Re_condition}). We thus have the following corollary.
\begin{cor}
If $\|\cdot\|_X$ is conjugation-invariant then $c_X \leq 2$. 
\end{cor}
We give a simple example of a norm that is not conjugation-invariant. It  also shows that the constant $c_X$ can be arbitrarily large. Let $r>0$ and define
$$\|x\|_X=r\,|x_1+ix_2|+|x_2|\qquad \text{for all }x=[x_1\;\; x_2]' \inn X\eqbd \C^2.$$
For $x=[1\;\; i]'$ we have $\|x\|_X=1$ and $\|\Re x\|_X=r$. Thus $c_X\geq r$.
\\

In order to determine the constant $c_X$ for the $L_p$-spaces we need
Jensen's inequality in the following form.
\begin{lemma}\label{lemma:Jensen}
({\bf Jensen's inequality)} Let $(\Omega, \cal B, \mu)$ be a measure space. Let $\kappa:[0,\infty)\to \R$ be a convex function. 
  Let $f,w\inn L^1_\R(\mu)$ be nonnegative functions  with 
$\int_\Omega w(t)\, d\mu_t=1$. Then 
$$\kappa\left(\int_\Omega f(t)\,w(t) \, d\mu_t\right)\leq
\int_\Omega \kappa(f(t))\,w(t) \, d\mu_t.$$
\end{lemma}
\begin{proof} For the special case $\mu(\Omega)=1$ and $w(t)\equiv 1$ the proof can be found in \cite[Theorem 3.3]{rudin}. To obtain the full statement define a measure $\tilde \mu$ on ${\cal B}$ by $\tilde \mu(B)\eqbd\int_Bw(t)\, d\mu_t$.
Then $\tilde \mu(\Omega)=1$ and 
$\int_\Omega f(t)\, d\tilde \mu_t= \int_\Omega f(t)\,w(t) \, d\mu_t$ for all
nonnegative measurable functions $f$
\cite[Theorem 1.29]{rudin}. Thus the result follows from the special case.
\end{proof}
To each measurable function $x:\Omega \to \C$ we associate a phase function
$\phi_x:\Omega \to [-\pi, \pi]$ defined by 
$$\phi_x(t)\eqbd \begin{cases} \pi, & \Im x(t)=0, \quad \Re x(t)\leq 0, \cr
2 \arctan ( {\Im x(t) \over \Re x(t)+\sqrt{(\Re x(t))^2+(\Im x(t))^2}}), & {\rm otherwise}. \end{cases} $$
Then $\phi_x$ is a measurable function, and, by elementary trigonometry,
 $x(t)=e^{i\phi_x(t)} |x(t)|$ for all $t \inn \Omega$.

The proposition below gives the constants $c_X$ for the $L_p$-spaces.
\begin{prop} Let $x\inn L_p^\C(\Omega, {\cal B}, \mu)\setminus \{ 0 \}$ where
$(\Omega, \cal B, \mu)$ is any measure space.  Then
\begin{eqnarray}\label{p_estimate}
\frac{\|\Re x\|_p+\|\Im x\|_p}{\|x\|_p}
\leq
\begin{cases}
\sqrt{2} &\text{if }1 \leq p \leq 2,\\
2^{1-1/p} &\text{if }2 \leq p < \infty,\\
2 &\text{if } p=\infty.
\end{cases} 
\end{eqnarray}
In the first case equality holds if $\Re x=\Im x$. In the second and the third
case equality holds if $\|\Re x\|_p=\|\Im x\|_p$ and $(\Re x(t)) (\Im x(t))=0$ for almost all $t\inn \Omega$.
\end{prop}
\begin{proof} 
We only show the estimate (\ref{p_estimate}). The proof of other statements
is left to the reader. The case $p=\infty$ is covered by case $(a)$ of
Proposition~\ref{prop:cx}.
Let $1\leq p \leq 2$. Then the function
$0\leq \xi \mapsto \xi^{2/p}$ is convex.
Let $w(t)=|x(t)|^p/\int_\Omega |x(\tau)|^p\, d\mu_\tau$.
 Then $\int_\Omega w(t)\, d\mu_t=1$. 
Applying Jensen's inequality, we obtain
\begin{eqnarray*}
\frac{\|\Re x\|_p^2 + \|\Im x\|_p^2}{\|x\|_p^2}
&=&
\left( \int_\Omega |\cos(\phi_x(t))|^p\, w(t)\, d\mu_t\right)^{2/p}
+ 
\left( \int_\Omega |\sin(\phi_x(t))|^p\, w(t)\, d\mu_t\right)^{2/p}
\\
&\leq&
 \int_\Omega |\cos(\phi_x(t))|^2\, w(t)\, d\mu_t
+
\int_\Omega |\sin(\phi_x(t))|^2\, w(t)\, d\mu_t
\\
&=&1.
\end{eqnarray*}
Thus,  
$\|\Re x\|_p + \|\Im x\|_p\leq \sqrt{2}\,\sqrt{\|\Re x\|_p^2 + \|\Im x\|_p^2}
\leq \sqrt{2}\,\|x\|_p.$

 Let $2 \leq p <\infty$. Then $(a^p+b^p)^{1/p}\leq (a^2+b^2)^{1/2}$ for
any $a,b\geq 0$. This implies $|\Re x(t)|^p+|\Im x(t)|^p\leq |x(t)|^p$ for
all $t \inn \Omega$. Consequently,
$\|\Re x\|_p^p+\|\Im x\|_p^p \leq \|x\|_p^p.$
By convexity, we also have
$\left(\frac{1}{2}\|\Re x\|_p+\frac{1}{2}\|\Im x\|_p\right)^p
\leq \frac{1}{2}\left(\|\Re x\|_p^p+\|\Im x\|_p^p\right).$
Combining the last two inequalities, we obtain
$$\|\Re x\|_p+\|\Im x\|_p \leq 2^{1-1/p}
\left(\|\Re x\|_p^p+\|\Im x\|_p^p\right)^{1/p} \leq 2^{1-1/p}\|x\|_p.$$
\quad \\[-1 cm]$\text{ }$
\end{proof}
\section{Cases of equality $\|A\|_{X,Y}^{\C}=\|A\|_{X,Y}^{\R}$}\label{sec:equal}
We first state the main result of this paper.
\begin{theo}\label{theo:main} 
Let $(\Omega_k,{\cal B}_k,\mu_k)$ be measure spaces, $k=1,2$. Let  
$A:L_p^\C(\mu_1)\to L_q^\C(\mu_2)$ 
be a linear operator that satisfies
$A(L_p^\R(\mu_1))\subseteq  L^q_\R(\mu_2)$. If $1\leq p\leq q \leq \infty$,
then  $\|A\|_{p, q}^{\R}=\|A\|_{p, q}^{\C}$.
\end{theo}
The proof is based on the following lemmas.
\begin{lemma}\label{lemma:x_phi_y}
Let $(X,X_\R,\|\cdot\|_X)$ and $(Y,Y_\R,\|\cdot\|_Y)$ be normed complexified
spaces. Let $A:X \to Y$ be a linear operator that satisfies 
$A(X_\R) \subseteq Y_\R$. Suppose that the following condition holds:
 For any $x \inn X$, $y \inn Y$ with $\|x\|_X=\|y\|_Y=1$ 
there exists a $\phi \inn [0,2\pi]$ such that
$\|\Re(e^{i\phi}x)\|_X\leq \|\Re(e^{i\phi}y)\|_Y$. Then 
$\|A\|_{X,Y}^\R=\|A\|_{X,Y}^\C.$
\end{lemma}
\begin{proof}
The case $A=0$ is trivial.
Let $x \inn X$ and suppose $Ax \not =0$.
By assumption, there exists $\phi \inn [0,2\pi]$ such that
\begin{eqnarray}\label{eq:phi_choice}
\left\|\Re\left(e^{i\phi}\frac{x}{\|x\|_X}\right)\right\|_X\leq 
\left\|\Re\left(e^{i\phi}\frac{Ax}{\|Ax\|_Y}\right)\right\|_Y
\end{eqnarray}
Let $\tilde x\eqbd\Re(e^{i\phi}x)\inn X_\R$. The condition 
$A(X_\R) \subseteq Y_\R$ yields that $\Re(e^{i\phi}Ax)=A\tilde x$.
Hence, it follows from (\ref{eq:phi_choice}) that
$\|Ax\|_Y/\|x\|_X\leq \|A \tilde x\|_Y/\|\tilde x\|_X.$ Thus 
$\|A\|_{X,Y}^\C\leq \|A\|_{X,Y}^\R$.
\end{proof}
\begin{lemma}\label{lemma:continuity}
Let $(X,X_\R,\|\cdot\|_X)$ be a normed complexified space. Suppose that
$\|\Re x\|_X \leq \|x\|_X$ for all $x \inn X$. Then for any $x \inn X$ 
the map $\phi \mapsto \|\Re(e^{i\phi}x)\|_X$, $\phi \inn \R$, is continuous.
\end{lemma}
\begin{proof} This follows from the inequalities
$$|\|\Re(e^{i\phi}x)\|_X-\|\Re(e^{i\phi_0}x)\|_X|
\leq \|\Re((e^{i\phi}-e^{i\phi_0})x)\|_X\leq 
\|(e^{i\phi}-e^{i\phi_0})x\|_X=|e^{i\phi}-e^{i\phi_0}|\, \|x\|_X,$$
which hold for all $\phi,\phi_0 \inn \R$
\end{proof}

\begin{lemma} \label{lemma:main_inequality}
Let $1 \leq p \leq q<\infty$. Then for any $x \inn L_p^\C(\mu)$ with $\|x\|_p=1$,
\begin{eqnarray}
\int_{0}^{2 \pi} \|\Re(e^{i\phi}x)\|_p^q \;\,  d\phi
&  \leq &   
  \int_{0}^{2 \pi} \!\! |\cos \phi|^q\, d\phi.
\end{eqnarray}
Equality holds if $p=q$. 
\end{lemma}
\begin{proof} First note that 
$\Re(e^{i\phi}x(t))=\cos(\phi_x(t)+\phi)|x(t)|$,
where $\phi_x$ is the phase function defined in Section~\ref{sec:bounds}.
Since $q\leq p$, the function $0\leq \xi \mapsto \xi^{q/p}$ is convex, 
and due to another assumption of the Lemma, $\int_\Omega |x(t)|^p\, d\mu_t =1$.
Hence, by Jensen's inequality as given in Lemma~\ref{lemma:Jensen},
for all $\phi \inn \R$,
\begin{eqnarray}
\nonumber
\|\Re(e^{i\phi}x)\|_p^q
&=&
\left(
\int_\Omega |\cos(\phi_x(t)+\phi)|^p\, |x(t)|^p
\, d\mu_t 
\right)^{q/p} \\[.2 cm]
\label{eq:leq1}
&\leq& 
\int_\Omega |\cos(\phi_x(t)+\phi)|^q\, |x(t)|^p
\, d\mu_t.
\end{eqnarray}
Combined with Fubini's Theorem, this yields 
\begin{eqnarray}\label{eq:leq2}
\int_{0}^{2\pi}\|\Re(e^{i\phi}x)\|_p^q
\, d\phi 
&\leq &
\int_{0}^{2\pi}
\int_\Omega |\cos(\phi_x(t)+\phi)|^q\, |x(t)|^p
\, d\mu_t\, d\phi\\[.2 cm]
&=&\nonumber
\int_\Omega
\int_{0}^{2\pi}
|\cos(\phi_x(t)+\phi)|^q\, |x(t)|^p
\,  d\phi \, d\mu_t\\[.2 cm]
&=&\nonumber
\int_\Omega
\left( \int_{0}^{2\pi}
|\cos(\phi_x(t)+\phi)|^q\, 
 d\phi\right)\, \,|x(t)|^p \, d\mu_t\\[.2 cm]
&=&\nonumber
\int_\Omega
\left(\int_{0}^{2\pi}
|\cos(\phi)|^q\, 
 d\phi\right)\, \, |x(t)|^p \, d\mu_t
\\[.2 cm]
&=&\nonumber
\int_{0}^{2\pi}|\cos(\phi)|^q\, d\phi.
\end{eqnarray}
If $p=q$, then equality holds in (\ref{eq:leq1}) and in (\ref{eq:leq2}).
\end{proof}

For the record, we next give an alternative proof of 
Lemma~\ref{lemma:main_inequality} for the finite-dimensional case.
It uses the H\"older norm 
$${\bf |\!\!|} f{\bf |\!\!|} 
\eqbd \left(\int_0^{2\pi} |f(\phi)|^{q/p}\, d\phi\right)^{p/q}, 
\qquad f \inn {\cal C}([0,2\pi], \R).$$

\begin{proofof} {\bf the finite-dimensional version of 
Lemma~\ref{lemma:main_inequality}:\/}
Write the  components of the vector $x\in \C^n$ in the form
$x_t= e^{i\phi_t}|x_t|$, $\phi_t \inn [0,2\pi]$, $t=1, \ldots, n$.
Then the $t$th component of the vector $\Re (e^{i\phi }x)$ is 
$|x_t|\cos(\phi_t+\phi)$. 
Let $f_t(\phi)\eqbd |x_t|^p\, |\!\cos(\phi_t+\phi)|^p.$ Then
\begin{eqnarray*}
{\bf |\!\!|}f_t{\bf |\!\!|} &=&
\left(\int_0^{2\pi} |f_t(\phi)|^{q/p}\, d\phi\right)^{p/q}\\
&=&
\left(\int_0^{2\pi} |x_t|^q\, |\!\cos(\phi_t+\phi)|^q\, d\phi\right)^{p/q}\\
&=&|x_t|^p\left(\int_0^{2\pi}|\cos(\phi)|^q\, d\phi \right)^{p/q}.
\end{eqnarray*}
Thus 
\begin{eqnarray}
\left(\int_0^{2\pi}\!\!\|\Re (e^{i\phi }x)\|_p^q\, d\phi
\right)^{p/q}   \nonumber
&=&
\left(\int_0^{2\pi}\!\!
\left(
\mbox{$\sum_{t=1}^n$}|x_t|^p\, |\cos(\phi_t+\phi)|^p
\right)
^{q/p}
\, d\phi
\right)^{p/q} \qquad\\[.2 cm] \nonumber
&=& 
\mbox{$
{\bf |\!\!|}\sum_{t=1}^nf_t{\bf |\!\!|}
$ }
\\[.2 cm] \label{eq:min2}
&\leq& 
\mbox{$
\sum_{k=t}^n{\bf |\!\!|}f_t{\bf |\!\!|}
$} 
\\ \nonumber
&=&
\|x\|_p^p\left(\int_0^{2\pi}|\cos(\phi)|^q\, d\phi \right)^{p/q}.
\end{eqnarray}
Thus
\begin{eqnarray}\label{eq:min3}
\int_0^{2\pi}\!\!\|\Re (e^{i\phi }x)\|_p^q\, d\phi
\leq \|x\|_p\, \int_0^{2\pi}|\cos(\phi)|^q\, d\phi.
\end{eqnarray}
If $p=q$,  then equality holds in 
(\ref{eq:min2}) and (\ref{eq:min3}). \end{proofof} 

We are now in a position to prove Theorem~\ref{theo:main} for $q < \infty$.
\vskip 0.2cm

\begin{proofof} {\bf Theorem~\ref{theo:main} for $q<\infty$:\/}
Let $1\leq p \leq q<\infty$. Let $x \inn L_\C^p(\mu)$, $y\inn L^q_\C(\mu)$ with 
$\|x\|_p=\|y\|_q=1$. Since the inequality $\|\Re(z)\|_p \leq \|z\|_p$
holds for all $z\inn L_p^\C(\mu)$ and all $p \inn [1, \infty]$, the function
$\phi \mapsto  \|\Re(e^{i\phi}y)\|_q^q- \|\Re(e^{i\phi}x)\|_p^q$ is continuous
by Lemma~\ref{lemma:continuity}. It follows from Lemma~\ref{lemma:main_inequality}
that
\begin{eqnarray*}
\int_0^{2\pi}(\|\Re(e^{i\phi}y)\|_q^q- \|\Re(e^{i\phi}x)\|_p^q)\, d\phi
&\geq & 0.
\end{eqnarray*}
Hence, the integrand is nonnegative for at least one $\phi_0 \inn [0,2\pi]$.
Thus $\|\Re(e^{i\phi_0}y)\|_q\geq \|\Re(e^{i\phi_0}x)\|_p$. Now,
Lemma~\ref{lemma:x_phi_y} yields the result.\end{proofof}

For $q=\infty$, the statement of Theorem \ref{theo:main} is covered by the 
following more general result.
\begin{theo}\label{theo:main2} 
Let $(\Omega,{\cal B},\mu)$ be a measure space. Let $(X,X_\R,\|\cdot\|_X)$
be a normed complexified space.  Suppose that $\|\Re x\|_X\leq \|x\|_X$
for all $x \inn X$. 
Let $A:X\to L_\infty^\C(\mu)$ 
be a linear operator that satisfies
$A(X_\R)\subseteq  L^\R_\infty(\mu)$.   
Then  $\|A\|_{X,\infty}^{\R}=\|A\|_{X, \infty}^{\C}$.
\end{theo}
For the proof we need yet another lemma.
\begin{lemma}\label{lemma:phi_selection}
 Let $x \inn L^\infty_\C(\mu)$. 
Then there is a $\phi_0 \inn [0,2\pi]$ 
with $\|\Re(e^{i\phi_0}x)\|_\infty=\|x\|_\infty$.
\end{lemma}
\begin{proof}
By Lemma \ref{lemma:continuity} the function 
$\phi \mapsto \|\Re(e^{i\phi}x)\|_\infty$, $\phi \inn [0,2\pi]$,
  is continuous. It therefore attains its maximum.  
Hence it is enough to show that to each $0<c<1$ there corresponds a 
$\phi_c \inn [0,2\pi]$ such that
$\|\Re(e^{i\phi_c}x)\|_\infty \geq c\, \|x\|_\infty$.
Let $M\eqbd 
\left\{\; t \inn \Omega\; |\;\; |x(t)|\geq \sqrt{c}\, \|x\|_\infty\;\right\}$. 
Note that $\mu(M)>0$, since $\sqrt{c}<1$. Choose an integer $n>1$ such that
$\cos(\theta)>\sqrt{c}$ for all $\theta\inn \R$ with 
$|\theta|\leq \frac{\pi}{n}$.
Since $M$ is the union of the sets
$$M_\ell=\left\{\;\; t \inn M\;\left|\;\; |\phi_x(t)-\frac{2\pi}{n}\ell |
\leq \frac{\pi}{n}\;\right.\right\}, \qquad \ell=1,\ldots ,n-1,$$
we have that $\mu(M_\ell)>0$ for at least one $\ell=:\ell_0$. Let 
$\phi_c=-\frac{2\pi}{n}\ell_0$. Then for any $t \inn M_{\ell_0}$,
we have $\cos(\phi_x(t)+\phi_c)>\sqrt{c}$, and hence
\begin{eqnarray*}
\Re(e^{i\phi_c}x(t))&=& \cos(\phi_x(t)+\phi_c)\, |x(t)|\geq c\|x\|_\infty.
\end{eqnarray*}
This yields $\|\Re(e^{i\phi_c}x)\|_\infty \geq c\, \|x\|_\infty$.
\end{proof}
\begin{proofof}{\bf Theorem \ref{theo:main2}:} 
Let $x \inn X$, $y \inn L_\infty^\C(\mu)$ with $\|x\|_X=\|y\|_\infty=1$.
By Lemma~\ref{lemma:phi_selection}, there is a $\phi_0$ with $\|\Re (e^{i \phi_0}y)\|_\infty=1$. Furthermore, we have
$\|\Re (e^{i \phi_0}x)\|_X\leq\|e^{i \phi_0}x\|_X=1$. Thus
$\|\Re (e^{i \phi_0}x)\|_\infty\leq \|\Re (e^{i \phi_0}y)\|_\infty.$
Now, the Theorem follows from Lemma \ref{lemma:x_phi_y}.
\end{proofof}

Theorem~\ref{theo:main2} can be proved more directly in a finite-dimensional
case. Then the norm on the target space, say $\C^m$, is nothing but a 
{\em weighted maximum norm,\/} i.e.,  a norm of the form
$$\|y\|_w=\max_{j\in \underline{m}} (w_j\, |y_j|),$$
where $w=[w_1 \ldots w_m]^T>0$ is a positive vector of weights.
The source space $\C^n$, however, does not have to be equipped with a
H\"older norm but, more generally, with any {\em absolute norm,\/} i.e., 
a norm that satisfies 
$$\|x\|=\|\, |x|\, \| \qquad \text{\rm for all } \quad x\in \C^n.$$
A n absolute vector norm  is necessarily {\em monotonic\/}~\cite{BSW}:
if $0\leq x\leq y$ componentwise, then $\|x\|\leq \|y\|$.
H\"older norms are obviously absolute.

\begin{prop}\label{prop:weights}
Let $\|\cdot\|_\alpha$ be an absolute norm on $\C^n$ and let $\|\cdot\|_w$ be a
weighted maximum norm on $\C^m$. 
Then $\|A\|_{\alpha, w}^{\R}=\|A\|_{\alpha, w}^{\C}$ for every 
$A =[a_{jk}]\inn \R^{m \times n}$. 
\end{prop} 
\begin{proof}
Let $x \inn \C^n$ and $j \inn \underline{m}$ be such that 
$\|Ax\|_w=w_j\left|\sum_{k \in \underline{n}}a_{jk}x_k\right|$.
Set $$ \widetilde{x}_k \eqbd \begin{cases} |x_k|, & a_{jk}\geq 0 \cr
-|x_k| & \text{\rm otherwise.} \end{cases} $$ 
Then we have for the vector $\widetilde{x}=
[\widetilde{x}_1 \ldots \widetilde{x}_n]'\inn \R^n$ that 
$\|\widetilde{x}\|_\alpha=\|\, |x|\,\|_\alpha=\| x\|_\alpha$ 
and  $\|A\widetilde{x}\|_w\geq
w_j\left|\sum_{j \in \underline{m}}a_{jk}\widetilde{x}_k\right|=w_j
\sum_{k\in \underline{n}}|a_{jk}x_k|\geq \|Ax\|_w.$ Thus
$\|A\widetilde{x}\|_{w}/\|\widetilde{x}\|_{\alpha}\geq \|Ax\|_{w}/\|x\|_{\alpha}.$
Thus $\|A\|_{\alpha, w}^{\R}\geq \|A\|_{\alpha, w}^{\C}$.
\end{proof}
Recall that the {\em dual norm\/}  $\|\cdot\|_\alpha^d$ associated with 
a given norm $\|\cdot\|_\alpha$ on $\C^n$  is defined by
\begin{eqnarray}\label{eq:dual}
\|a\|_\alpha^d=\max_{x \in \C^n\setminus \{0\}}
\frac{|a' x|}{\|x\|_\alpha}\qquad \text{for }a \in \C^n.
\end{eqnarray}
By specializing Proposition~\ref{prop:weights} to the case $m=1$, we obtain
the following well-known fact~\cite{BSW}.
\begin{cor}
Let $\|\cdot\|_\alpha$ be an absolute norm on $\C^n$.
If $a \inn \R^n$, then the maximum in (\ref{eq:dual}) is attained for a 
real vector $x$.
\end{cor} 

The next result is a variant of Theorem \ref{theo:main} for
inner product spaces.
\begin{theo}
Let $(X,X_\R,\|\cdot\|_X)$, $(Y,Y_\R,\|\cdot\|_Y)$ be  normed complexified space with
\begin{eqnarray*}
\|x\|_X^2&=&\langle \Re x ,\Re x\rangle +  \langle \Im x ,\Im x\rangle
\qquad \text{\rm for  } x \inn X,\\
\|y\|_Y^2&=&( \Re x ,\Re x) +  (\Im y ,\Im y)
\qquad \text{\rm for  } y \inn Y,
\end{eqnarray*}
where $\langle \cdot, \cdot \rangle$ and $(\cdot , \cdot )$ are inner products
on on $X_\R$ and $Y_\R$ respectively.
Then $\|A\|_{X,Y}^\R=\|A\|_{X,Y}^\C$ for any linear operator $A:X \to Y$
satisfying $A(X_\R)\subseteq  Y_\R$.
\end{theo}
\begin{proof} A straightforward computation yields that for any $x\inn X$,
$$\int_{0}^{2 \pi} \|\Re(e^{i\phi}x)\|_X^2 \;\,  d\phi=\left(\int_{0}^{2 \pi}
\!\!\!\cos^2(\phi)\, d\phi\right)\|x\|_X^2=\pi \, \|x\|_X^2.$$
The same relation holds for $y\inn Y$ and $\|\cdot\|_Y$. Thus, if $\|x\|_X=\|y\|_Y=1$,
$$\int_{0}^{2 \pi} 
(\|\Re(e^{i\phi}y)\|_Y^2-\|\Re(e^{i\phi}x)\|_X^2)\, d\phi=0.$$
Hence, by continuity $\|\Re(e^{i\phi}x)\|_X\leq \|\Re(e^{i\phi}y)\|_Y$
 for some $\phi \inn [0, 2\pi]$, and Lemma \ref{lemma:x_phi_y} applies.
\end{proof}

We close the section by showing that the real and complex norms of a nonnegative
linear operator coincide whenever the source $X$ and the target $Y$ are complexified
function spaces where $\|\cdot \|_X$ is absolute and $\|\cdot \|_Y$ is absolute and 
monotonic. This means that $X_\R$ and $Y_\R$ are spaces of real-valued functions 
and, for  each $x\in X$ ($y\in Y$), the absolute 
value function  $|x|\eqbd \sqrt{(\Re x)^2+(\Im x)^2}$ 
($|y|\eqbd \sqrt{(\Re y)^2+(\Im y)^2}$) is an element of the space 
$X_\R$ ($Y_\R$) and $\|x\|_X=\| \,|x|\, \|_X$ ($\|y\|_Y=\|\, |y|\, \|_Y$). 
Moreover, if $f,g\in Y_\R$ and $0\leq f \leq g$ pointwise, then 
$\| f\|_Y \leq \| g\|_Y$. 
An operator from $X$ to $Y$ is {\em nonnegative\/} if it is real and it
maps nonnegative functions from $X_\R$ to nonnegative functions from $Y_\R$.

\begin{theo}\label{theo:nonnegative} Let $(X,X_\R,\|\cdot\|_X)$ and
 $(Y,Y_\R,\|\cdot\|_Y)$ be complexified function spaces, let 
$\|\cdot \|_X$ be absolute and $\| \cdot \|_Y$ be absolute and 
monotonic, and let $A : X \to Y$ be a nonnegative linear operator. Then
$\|A\|^\C_{X,Y}=\|A\|^\R_{X,Y}$. \end{theo}

\begin{proof} First note that the absolute value function is defined by the property 
$$|f|=\inf \{g\geq 0 : g\geq \Re (zf) \; {\rm whenever} \; z\in \C, \; |z|=1\}.$$
Also observe that $A$, being a real operator,  commutes with taking 
the real value of a function. Since $A$ is nonnegative nad linear, 
$Af\geq Ag$ whenever $f\geq g$, hence
$$\Re(zAf)=A\Re(zf)\leq A |f| \quad {\rm whenever} \quad z\in \C,\; |z|=1.$$
This implies that $|Af|\leq A|f|$. By the absoluteness and monotonicity of~$\|\cdot \|_Y$,
we have $\|Af\|_Y=\| \,|Af|\, \|_Y\leq \|\,A|f|\,\|_Y$, whereas $\|f \|_X=\|\, |f| \,\|_X$ 
since $\|\cdot\|_X$ is absolute. The function $f\in X$ was arbitrary, so
the norm $ \|A\|^\C_{X,Y}$ is in fact equal to 
$$ \sup_{0\neq f\in X_R, \; f=|f|} { \| Af\|_Y\over \|f\|_X}  $$
and therefore to  $ \|A\|^\R_{X,Y}$. \end{proof}
 
Alternatively, in the finite-dimensional case, the equality of real and complex 
norms  can be seen as follows. We denote by $|A|$ the matrix of absolute values of 
$A \inn \C^{m\times n}$.

\begin{prop}
Let $\|\cdot\|_\alpha$, $\|\cdot\|_\beta$  be absolute norms on $\C^n$ 
and $\C^m$ respectively.  Then 
$\|A\|_{\alpha, \beta}^\C \leq \|\; |A|\; \|_{\alpha, \beta}^\R$
for all $A \inn \C^{m \times n}$.
\end{prop}
\begin{proof}
The  monotonicity property of absolute norms yields  that for any
 $x \inn \C^n$,
\begin{eqnarray*}
\| Ax \|_\beta &=& \| \; |Ax|\;  \|_\beta\\
&\leq& \| \; |A|\,|x|\;  \|_\beta\\
&\leq& \| \; |A|\;\|_{\alpha,\beta}^\R \; \|\; |x|\;\|_\alpha\\
&=&\| \; |A|\;\|_{\alpha,\beta}^\R \;\|x\|_\alpha.
\end{eqnarray*}
The result follows. 
\end{proof}
\begin{cor} If a matrix $A \in \R^{m \times n}$ is elementwise nonnegative and the underlying 
norms $\|\cdot \|_\alpha$ on $\C^n$ and $\| \cdot \|_\beta$ on $\C^m$ are absolute, then 
$\|A\|_{\alpha, \beta}^\C = \|A \|_{\alpha, \beta}^\R.$
\end{cor}

\section{Cases of inequality $\|A\|_{X,Y}^{\C}>\|A\|_{X,Y}^{\R}$}\label{sec:unequal}
We now show that the main result, Theorem \ref{theo:main},  is sharp
already in the finite-dimensional case. In other words, for any $p>q\geq 1$,
there exists a real matrix $A$ such that
$$  \|A\|^\C_{p,q}>\|A\|^\R_{p,q}. $$
We begin with the case $p>q$, $q<2$.

\begin{prop} Let $p>q$, $q<2$, and let
$$ A\eqbd \left[ \begin{array}{rrr} 1 & 1 & 0  \\ 1 & -1 & 0  \\
0 & 1 & 1  \\ 0 & 1 & -1  \\ 1 & 0 & 1  \\ 1 & 0 & -1  \end{array} \right].$$
Then $ \|A\|^\C_{p,q}>\|A\|^\R_{p,q}.$
\end{prop}

\begin{proof} First let us show that the value of the real 
$(p,q)$-norm of the matrix $A$ is {\em not\/} attained at the
vector with just one nonzero component. Due to the symmetry of the 
entries of $A$, it is enough to argue about the first unit vector 
$v\eqbd[1,0,0]'$. Consider its small real perturbation $v(\varepsilon)\eqbd
[1,\varepsilon,0]'$ and compare corresponding ratios of norms:
\begin{eqnarray*}
{\| Av\|^q_q  \over   \|v\|_p^q} & = & 4 \\
{\| Av(\varepsilon)\|^q_q  \over  
\|v(\varepsilon)\|_p^q} & = & 4+2\varepsilon^q-4{q\over
 p}\varepsilon^p +q(q-1)\varepsilon^2+o(\varepsilon^2).  \\
\end{eqnarray*}
The latter expression is strictly greater than $4$ for small $\varepsilon$,
since $2\varepsilon^q$ is then the smallest order term. 
Thus, the real $(p,q)$-norm of $A$ is attained at a vector
with at least two nonzero components, say $v_{\min{}}\eqbd[x,y,z]'$.
Again, due to the form of the matrix $A$, the components can
be assumed to be all nonnegative and ordered so that $x\geq y\geq z\geq 0$,
$y>0$.

Now, consider the vector $w\eqbd [ix,y,z]'$. Since the function 
$f(x)\eqbd x^{q/2}$ is strictly concave on the nonnegative real
axis, we have
\begin{eqnarray*}
(x^2+y^2)^{q/2} &>&{ ((x+y)^2)^{q/2}+ 
(x-y)^2)^{q/2} \over 2}, \\
(x^2+z^2)^{q/2}& \geq &{ ((x+z)^2)^{q/2}+ 
(x-z)^2)^{q/2} \over 2}, \\
\end{eqnarray*}
hence $\|Av_{\min{}}\|_q^q<\| Aw\|_q^q$, whereas
$\|v_{\min{}}\|_p=\|w\|_p$. Hence, the complex $(p,q)$-norm 
of $A$ is strictly bigger than its real norm.
\end{proof}

The case $p>q\geq 2$ reduces to the case we just considered 
due to duality. We state this formally as a lemma.

\begin{lemma} Suppose that a real matrix $A$ satisfies 
$ \|A\|^\C_{p,q}>\|A\|^\R_{p,q}. $
Then its transpose $A'$ satisfies
$ \|A'\|^\C_{q',p'}>\|A'\|^\R_{q',p'}, $
where $p'\eqbd{p/ (p-1)}$, $q'\eqbd{q/ (q-1)}$.
\end{lemma}

\begin{proof} Follows directly from the fact
$$ \|A\|^\F_{p,q}=\|A'\|^\F_{q',p'}, $$
which holds  for  both $\F=\R$ and $\F=\C$.
\end{proof}

Since $p'<2$ whenever $p>2$, this observation   
enables us to produce counterexamples for the case $p>q\geq 2$
out of couterexamples for the previous case. 

\begin{cor} Let $p>q\geq 2$, and let
$$ A\eqbd \left[ \begin{array}{rrrrrr} 1 & 1 & 1 & 1 & 0 & 0 \\ 
1 & -1 & 0 & 0 & 1 & 1 \\ 0 & 0 &  1 & -1  & 1 & -1   \end{array} \right].$$
Then $ \|A\|^\C_{p,q}>\|A\|^\R_{p,q}.$
\end{cor}

This finishes our proof that the condition $p\leq q$ 
in the main theorem of this paper cannot be relaxed.
\section{Norm extensions}\label{sec:extend}
In this section we provide a class of norm extensions from a real
vector space to its complexification which preserve operator norms.

Let $\nu$ be a norm on ${\cal C}([0,2\pi],\R)$,
the set of continuous real valued functions
on the interval $[0,2\pi]$. Then  $\nu$ is said to be monotone if
$0\leq f(t)\leq g(t)$  for all $t \inn [0,2\pi]$ implies that
$\nu(f) \leq \nu(g)$. The norm $\nu$ is called shift-invariant 
if for all $\psi \in \R$, $\nu(f)=\nu(f_\psi)$, where
$f_\psi(\phi)\eqbd f(\, (\phi+\psi) \,{\rm mod}\, 2\pi)$. 
For instance, the $L_p$-norms are monotone and shift-invariant.

Let $X_\R$ be a vector space over $\R$ endowed with a norm $\|\cdot\|_{X_\R}$.
Let $X$ be the algebraic complexification of $X_\R$. It is then easily
seen that for any $x \inn X$ the function 
$\phi \mapsto \|\Re (e^{i\phi} x)\|_{X_\R}$, $\phi \inn [0, 2\pi]$, 
is continuous. 
We denote these functions by $\|\Re (e^{i\cdot} x)\|_{X_\R}$.
Let $\nu$ be a monotone and shift-invariant norm on ${\cal C}([0,2\pi],\R)$.
For $x \inn X$ define
\begin{eqnarray}\label{eq:xnorm}
\|x\|_X\eqbd \nu(\|\Re (e^{i\cdot} x)\|_{X_\R}).
\end{eqnarray}
\begin{prop}
The function $\|\cdot\|_X:X \to \R$ is a norm on $X$. 
Suppose $\nu$ is normalized so that $\nu(\, \,|\!\cos(\cdot)|\,\,)=1$.
Then $\|x\|_X=\|x\|_{X_\R}$ for all $x\inn X_\R$. In other words
 $\|\cdot\|_X$ is an extension of $\|\cdot\|_{X_\R}$.
\end{prop}
\begin{proof}
The triangle inequality for $\|\cdot\|_X$ follows from the triangle
inequality for $\|\cdot\|_{X_\R}$ and the monotonicity of $\nu$. 
The identity $\|\lambda \, x\|_X=|\lambda|\, \|x\|_X$ for all 
$\lambda \inn \C$ is a consquence of the shift invariance of $\nu$. 
The rest is obvious.
\end{proof}
In the proposition below  $Y_\R$ is a second vector space over $\R$
endowed with a norm $\|\cdot\|_{Y_\R}$ and complexification $Y$. 
The norm $\|\cdot\|_Y$ is defined by the same $\nu$, i.e.
\begin{eqnarray}\label{eq:ynorm}
 \|y\|_Y\eqbd \nu(\|\Re (e^{i\cdot} y)\|_{Y_\R})
\end{eqnarray}
\begin{prop} 
Let the norms on $X,Y$ be defined as in (\ref{eq:xnorm}), (\ref{eq:ynorm}),
where $\nu$ is monotone and shift-invariant.
Then $\|A\|_{X,Y}^\C=\|A\|_{X,Y}^\R$
for all linear maps $A:X \to Y$ satisfying $A(X_\R)\subseteq Y_\R$.
\end{prop}  
\begin{proof} We have,
\begin{eqnarray*}
\|Ax\|_Y &=& \nu(\|\Re(e^{i \cdot}Ax)\|_{Y_\R})\\
&= &\nu(\|A\Re(e^{i \cdot}x)\|_{Y_\R})\\
&\leq & \nu(\|A\|_{X,Y}^\R\, \|\Re(e^{i \cdot}x)\|_{X_\R})\\
&= &\|A\|_{X,Y}^\R \, \nu(\|\Re(e^{i \cdot}x)\|_{X_\R})\\
&=&\|A\|_{X,Y}^\R\, \|x\|_X.
\end{eqnarray*}
\quad \\[-1 cm]$\text{ }$
\end{proof}


\section*{Acknowledgements}

Part of this work was done while the first author was visiting 
Technische Universit\"at Berlin. She is grateful to her host,
Volker Mehrmann, for his hospitality. The authors also thank Hans 
Schneider for his comments on an earlier version of this paper.


\section*{Postscriptum}

After this paper was complete, the authors were informed by  Lech Maligranda 
that the main results of this paper had been published in~\cite{GaschMali}
and~\cite{Mali}.


\end{document}